\documentclass[12pt,amstex]{article}
\usepackage[usenames]{color}

\usepackage{amscd}
\usepackage{amsmath}
\usepackage{amssymb}
\usepackage{graphicx}

\newtheorem {theo} {\bf Theorem} [section]
\newtheorem {prop} [theo] {\bf Proposition}
\newtheorem {cory} [theo] {\bf Corollary}
\newtheorem {lem} [theo] {\bf Lemma}
\newtheorem {defn} [theo] {\bf Definition}

\newtheorem {rem} [theo] {\bf Remark}
\newcommand{\QED}{\hfill \CaixaPreta \vspace{6mm}}
\def\CaixaPreta{\vrule Depth0pt height6pt width6pt}

\newcommand{\qed}{\nopagebreak\hfill{\vrule width6pt height6pt depth0pt}}
\newcommand{\be}{\begin{eqnarray}}
\newcommand{\ee}{\end{eqnarray}}
\newcommand{\benn}{\begin{eqnarray*}}
\newcommand{\eenn}{\end{eqnarray*}}
\newcommand{\bse}{\begin{equation}}
\newcommand{\ese}{\end{equation}}
\newcommand{\bsenn}{\begin{displaymath}}
\newcommand{\esenn}{\end{displaymath}}
\newcommand{\logand}{\;\;{\rm and }\;\;}

\newcommand{\where}{\;\;{\rm where }\;\;}

\newcommand{\R}{\mathbb{R}}

\usepackage{graphicx}

\begin{document}

\title{Geodesics on Regular Constant Distance Surfaces}
\author{J. J. P. Veerman\thanks{Weizmann Institute, Rehovot, Israel.}
	\thanks{Maseeh Dept. of Math. and Stat., Portland State Univ., Portland, OR, USA; e-mail: veerman@pdx.edu.}\\
}\maketitle

\begin{abstract} Suppose that the surfaces $K_0$ and $K_r$ are the boundaries of two convex,
complete, connected $C^2$ bodies in $\R^3$. Assume further that the (Euclidean) distance between
any point $x$ in $K_r$ and $K_0$ is always $r$ ($r>0$).
For $x$ in $K_r$, let $\Pi(x)$ denote the nearest point to $x$ in $K_0$. We show that the projection
$\Pi$ preserves geodesics in these surfaces if and only if both surfaces are concentric spheres or
co-axial round
cylinders. This is optimal in the sense that the main step to establish this result is false for
$C^{1,1}$ surfaces. Finally, we give a non-trivial example of a geodesic preserving projection of two
$C^2$ \emph{non}-constant distance surfaces. The question whether for any $C^2$ convex surface $S_0$,
there is a surface $S$ whose projection to $S_0$ preserves geodesics is open.
\end{abstract}


\vskip 0.2in
\begin{centering}\section{Introduction}
	\label{chap:definitions}\end{centering}
\setcounter{figure}{0} \setcounter{equation}{0}

Suppose $\gamma(t)$ is a trajectory of an object in $\R^3$ outside a convex body.
In this paper, $\Pi(\gamma(t))$ is called the \emph{projection} of $\gamma(t)$. In many applications
it is important to track the point $\Pi(\gamma(t))$ on the surface of the body nearest to the moving
object\footnote{A case in point is the event on September 26, 2022, when an unmanned spacecraft hit the
asteroid Didymos on purpose \cite{naidu}, thereby changing the orbit of the asteroid.
Clearly, the change in orbit of the asteroid is related to the locus, angle, and speed of the missile
at the time of impact. The asteroid itself is in good approximation a convex set, but far from round \cite{naidu}.}.
In \cite{navigate}, a method to compute and track the projection was considered. Instead,
here we consider the question whether this projection can take geodesics to (reparametrized) geodesics.

Before describing the main result, we give some general background about this problem. A diffeomorphism
$\phi:S_1\rightarrow S_2$ between (sub) manifolds is called a \emph{geodesic mapping} if it carries
geodesics to geodesics. We restrict our discussion to surfaces in $\R^3$. It is well-known that if
$S_1$ has constant Gaussian curvature, then there is a geodesic mapping from $S_1$ to the plane.
Vice versa, Beltrami's theorem says that if $S_1$ admits a (local) geodesic mapping to the plane
near every point in $S_1$, then $S_1$ has constant Gaussian curvature (\cite{docarmo},
Section 4.6, exercises 12 and 13). There is a fairly large body of literature on geodesic mappings,
\cite{mikes1,mikes2} and the references therein. Our own interest here is to find out whether
\emph{projections} from one surface to another can be \emph{geodesic mappings}.

Our main result concerns projections to the convex set from a surface whose distance to the convex set
is exactly $r$ (a constant). We call such a surface a surface of constant distance (the word
`equidistant' is already in use for a slightly different concept \cite{wilker}). Very little has been
written about sets of constant distance (but see \cite{brown, misiure}). What we aim to show here
is essentially a rigidity result in $\R^3$: a constant distance surface whose projection takes
geodesics to geodesics must be a sphere or a cylinder. We proceed with the details.

We imagine a $C^2$ convex body in $\R^3$ whose boundary we denote by $K_0$. Let $p$ be any point
in the surface. By applying an isometry, we may assume that $p$ is located at the origin of $\R^3$
and that the tangent plane to $K_0$ at $p$ is given by $z=0$. Thus the coordinate patch near the
origin can be written as
\bse
K_0(x_1,x_2)=\left(x_1,x_2,-\frac 12(a_1x_1^2+a_2x_2^2) - h(x_1,x_2)\right)\,,
\label{eq:coord-patch1}
\ese
where the $a_i$ are the \emph{principal curvatures} and $h$ is twice continuously differentiable
with $h(0,0)$ is zero and the same holds for all first and second derivatives. By convexity, the
principal curvatures $a_i$ are non-negative.

Because of the smoothness and the convexity, we can smoothly coordinatize the space $\Omega$ surrounding
the convex body by using these coordinate patches as follows \cite{navigate}:
\bse
S(x_1,x_2,r)=K_0(x_1,x_2)+r\hat n(x_1,x_2) \,.
\label{eq:milnorcoords}
\ese
where $\hat n$ is the unit normal to $K_0$. These 3-dimensional coordinate patches form a differentiable
atlas of $\Omega$. Denote by $\Pi:\Omega\rightarrow K_0$ the orthogonal `projection'
from $\Omega$ onto $K_0$, defined as follows \cite{navigate}: $\gamma:=\Pi(z)$ is the unique point on $K_0$
nearest to $z\in \Omega$. Clearly, the inverse of $\Pi$ at a point
$\gamma$ of $K_0$ consists of a ray normal to $K_0$ at $\gamma$.
\bsenn
\Pi^{-1}(\gamma)= \cup_{r>0}\{\gamma+r\hat n(\gamma)\} \,,
\esenn
where $\hat n$ is the unit normal at $\gamma$ pointing outwards.

\begin{figure}[!ht]
\centering
\includegraphics[width=3.4in]{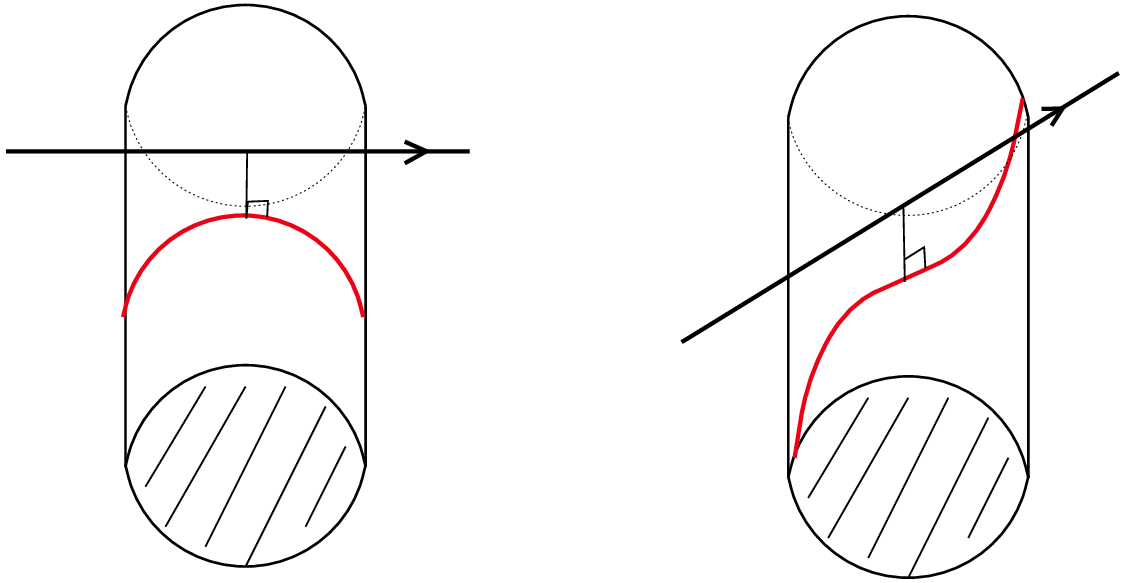}
\caption{\emph{Left, the projection (red) of straight line orthogonal to the axis of symmetry of a
solid cylinder. Right, the projection of a line at an arbitrary angle with the axis of symmetry of the
cylinder. The former is a geodesic, the latter clearly not.}}
\label{fig:projection-to-cyl}
\end{figure}

The simple example of Figure \ref{fig:projection-to-cyl} shows that the projection of a straight line in
$\R^3$ does not usually result in a geodesic in $K_0$. The question arises when is it that geodesics
\emph{do} project to geodesics? In this paper, any non-singular reparametrization (i.e.
with non-zero, possibly variable, speed) of a unit speed geodesic will also be called a geodesic.

\begin{defn} i) Given two closed $C^2$ surfaces $S_1$ and $S_2$ in $\R^3$. The \emph{projection}
$\Pi :S_1\rightarrow S_2$ is defined as follows. For $x\in S_1$,
\bsenn
\Pi(x):=\{y\in  S_2\,:\, y \textrm{ minimizes the Euclidean distance } d(x,y) \} \,.
\esenn
ii) Two surfaces are called \emph{regular constant distance surfaces} if the Euclidean distance from any
point $x$ in $S_1$ to $S_2$ equals $r$ (fixed), and the nearest point on $S_2$ is always unique.
\label{def:proj}
\end{defn}

\vskip 0.0in\noindent
It is a curious fact that in general $\Pi :S_1\rightarrow S_2$ and $\Pi' :S_2\rightarrow S_1$ are
\emph{not} inverses of one another. However, if the $S_i$ are at least $C^1$ and regular constant
distance, then $\Pi$ and $\Pi'$ \emph{are} inverses. This is the content of Proposition \ref{prop:projinverses}.
In the remainder of this paper, we deal with this case (except where mentioned otherwise).

Let $K_r$ denote the surface that has distance $r$ to $K_0$, or
\bsenn
K_r:=\{S(x_1,x_2,r)\,: \, r>0 \textrm{ fixed} \}\,.
\esenn
We are interested in determining when the projection
$\Pi:K_r\rightarrow K_0$ between these surfaces have the property that they send geodesics to
(reparametrizations of) geodesics. We call this property \emph{preservation of geodesics} and $\Pi$
a \emph{geodesic mapping}. The proof of the following result takes up most of this paper.

\begin{theo} Let $K_0$ be  $C^2$ surface patch given by \eqref{eq:coord-patch1} with $a_1\geq 0 $
and $a_2\geq 0$ and fix $r>0$. Then the projection $\Pi:K_r\rightarrow K_0$ does \emph{not} preserve
geodesics, unless (in that patch) (i) the Gaussian curvature is zero (i.e. $a_1a_2=0$)
or (ii) the patch consists of umbilic points (i.e. $a_1=a_2$)
\label{thm:funny geodesic2}
\end{theo}

A moment's reflection, will tell us that in $\R^3$, projections between concentric spheres or
between co-axial round cylinders \emph{do} preserve geodesics. The interesting question is, are
there any others? Here is a (to the author) surprising corollary of Theorem \ref{thm:funny geodesic2}.

\begin{cory} Let $K_0$ and $K_r$ be regular constant distance, complete, convex, connected, $C^2$ surfaces
in $\R^3$ at a distance $r>0$. The projection from $K_r$ to $K_0$ preserves geodesics if and only if
both are either spheres or (infinite) round cylinders.
\label{cory:funnygeodesic}
\end{cory}

\vskip 0.1in\noindent
{\bf Remark.} In this context, a (generalized) cylinder $C$ is a set of points such that for every point
$p\in C$ there is a unique line $\ell(p)$ in $C$ and any two such lines are either the same or parallel.
A `perfect' or `round' cylinder is a cylinder that rotationally symmetric around its axis. In particular, its
principal curvatures are constant.

\vskip 0.1in\noindent
{\bf Remark.} In view of Proposition \ref{prop:projinverses}, $\Pi:K_r\rightarrow K_0$ and
$\Pi':K_0\rightarrow K_r$ are inverses. So $\Pi$ preserves geodesics if and only if $\Pi'$
preserves geodesics.

\vskip 0.1in\noindent
{\bf Proof of Corollary \ref{cory:funnygeodesic}.} It is clear that if $K_0$ and $K_r$ both are either
spheres or (infinite) round cylinders, then the geodesics are preserved.

Vice versa, if the projection preserves geodesics, then by Theorem \ref{thm:funny geodesic2}, every
$C^2$ surface patch is either a piece of a sphere or piece of a cylinder. The two cannot occur
in the same $C^2$ patch, because at any `intermediate' point, (i) or (ii) in that theorem
will be violated, and then geodesics will not be preserved. Thus all of $K_0$ must satisfy either
(i) or (ii).

It is well-known that a $C^2$ complete surface whose principal curvatures are the identical
(or \emph{umbilic} surface) must be a part of a sphere (\cite{docarmo}, Section 3.2).
Similarly (\cite{docarmo}, section 5.8), a complete surface with Gaussian curvature zero,
must be a generalized cylinder. Finally, Proposition \ref{prop:round}
implies that if $K_0$ and $K_r$ are cylinders and the projection preserves geodesics, then
they must be round cylinders.
\QED

\vskip -0.1in\noindent
{\bf Remark.} Interestingly, this corollary is clearly false in $\R^2$. For instance, if $K_0$ is an ellipse
in $\R^2$ and $K_r$ a circle that contains it, the projection $K_r\rightarrow K_0$ is surjective.
On the other hand, in dimension 4 or higher, nothing appears to be known.

\vskip 0.1in
We furthermore prove that Corollary \ref{cory:funnygeodesic} is optimal in the sense that if we drop
$C^2$ in favor of $C^{1,1}$, that is: once continuously differentiable with a Lipschitz derivatives,
then the result does not hold.

\vskip-0.1in\noindent
\begin{theo} There exist regular constant distance, complete, convex, $C^{1,1}$ surfaces $K_0$ and $K_r$
in $\R^3$ with the property that (wherever the surfaces are $C^2$) either (i) $a_1a_2=0$ or (ii)
$a_1=a_2$ holds, but the projection from $K_0$ to $K_r$ does not preserve geodesics.
\label{thm:funny geodesic3}
\end{theo}

\vskip -0.1in\noindent
{\bf Proof.} The result follows directly from Proposition \ref{prop:funny geodesic3}.
\QED

\vskip-0.1in\noindent
{\bf Remark.} In \cite{Akmal} (see also \cite{navigate}), a related, but more complicated, counter-example
was constructed which carries over to cylinders in $\R^3$. It says that here is a convex $C^{1,1}$
cylinder such that the projection $\Pi$ onto this cylinder does not have a derivative.

\vskip0.1in
Finally, we are interested in the question whether, given the boundary $S_0$ of a convex body,
there is \emph{any} surface $S$ outside it, whose projection onto $S_0$ preserves geodesics.
For cylinders in $\R^3$, the answer is
affirmative, as we show in Section \ref{chap:preservegeod}. In fact, in that case, the space
outside $S_0$ can be foliated by surfaces $S_k$, $k\geq 0$ so that each projection $\Pi_k:S_k\rightarrow S_0$ preserves geodesics. However, as we will show, these surfaces $S_k$ generally are not convex.

\vskip 0.1in\noindent
{\bf Remark.} For general $C^2$ convex bodies, even in $\R^3$, it is unknown at the time of this writing
whether the space outside
them can be foliated by surfaces $S_k$ so that each projection $\Pi_k:S_k\rightarrow S_0$ preserves geodesics.

\vskip 0.2in
\begin{centering}\section{Preliminaries}
	\label{chap:outline}\end{centering}
\setcounter{figure}{0} \setcounter{equation}{0}

We first prove that the projections between two regular constant distant surfaces (see Definition
\ref{def:proj}) are inverses of one another. Then we discuss
the strategy to prove Theorem \ref{thm:funny geodesic2}.

\begin{prop} Let $S_1$ and $S_2$ be $C^1$ surfaces in $\R^3$ such that the Euclidean distance from any
point $x$ in $S_1$ to $S_2$ equals $r$ (fixed), and the nearest point on $S_2$ is always unique.
Then the projections $\Pi :S_1\rightarrow S_2$ and $\Pi' :S_2\rightarrow S_1$ are inverses of one
another.
\label{prop:projinverses}
\end{prop}

\begin{figure}[!ht]
\centering
\includegraphics[width=2.2in]{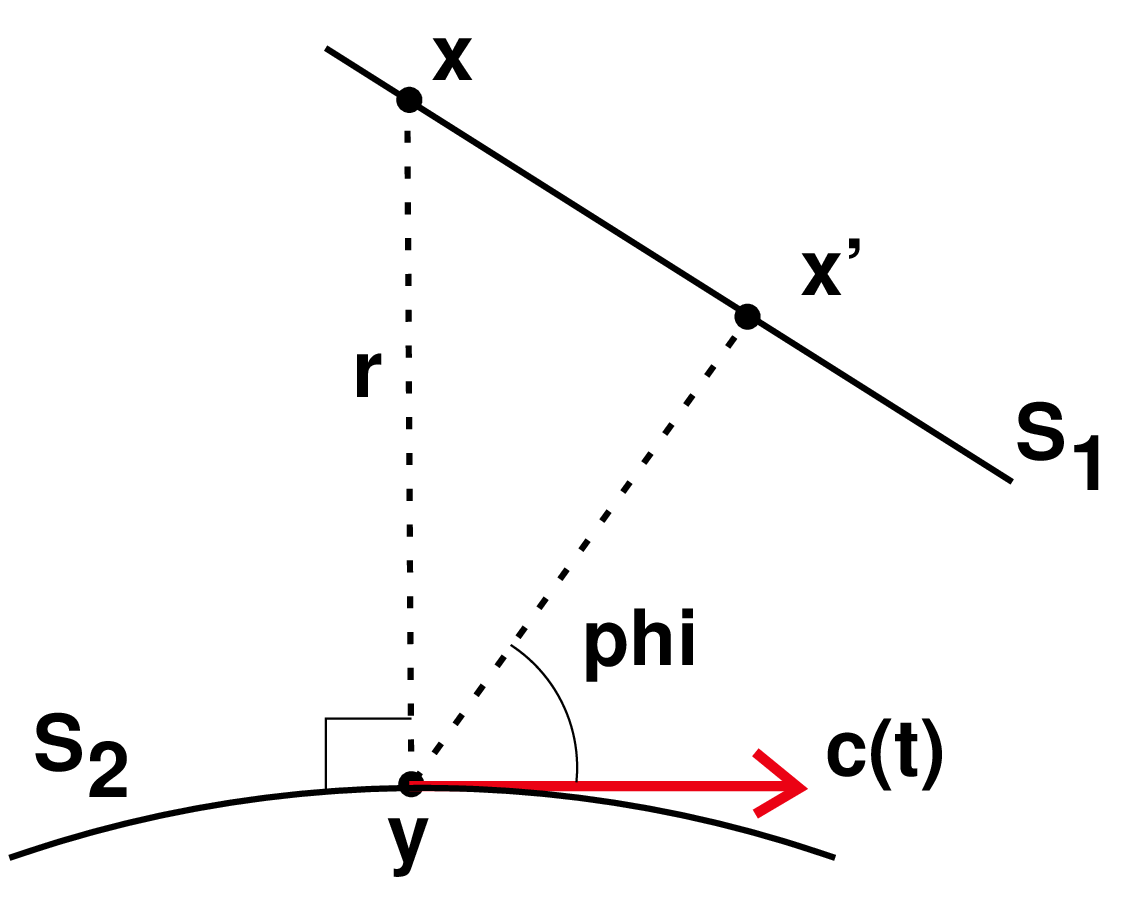}
\caption{\emph{This figure illustrates that $\Pi :S_1\rightarrow S_2$ and $\Pi' :S_2\rightarrow S_1$
are not generally inverses of one another. Traveling from $y$ along $S_2$ in the direction of $c(t)$
will (initially) decrease the distance to $x'$.}}
\label{fig:projinverses}
\end{figure}

\vskip 0.0in\noindent
{\bf Proof.} Consider $\Pi :S_1\rightarrow S_2$ and $\Pi' :S_2\rightarrow S_1$ and suppose $\Pi(x)=y$
(see Figure \ref{fig:projinverses}). Suppose there is $x'$ in $S_2$ not equal to $x$ such that
$x'\in \Pi'(y)$. Denote the Euclidean distance by $d(x,y)$. Now
\benn
x'\in \Pi'(y) \quad &\Longrightarrow& \quad d(y,x')\leq d(y,x)=r \\
d(x',S_2)=r \quad &\Longrightarrow& \quad d(x',y)\geq r
\eenn
So $d(y,x')=r$.

Consider the plane $P$ through $x$, $x'$, and $y$, and parametrize $S_2$ by the arclength $t$
and let the geodesic $c(t)$ be the tangent to $S_2(t)$ as drawn in Figure \ref{fig:projinverses}.
Then, by differentiability of $S_2$,
\bsenn
\lim_{t\searrow 0^+} \frac{d(S_2(t),x')-d(S_2(0),x')}{t}=
\lim_{t\searrow 0^+} \frac{d(c(t),x')-d(c(0),x')}{t}=-\cos \phi \,.
\esenn
The last equality is a special case\footnote{In this simple case, it can also be derived easily from an
explicit computation.} of Theorem 4.3 in \cite{Fox}. Thus for some positive $t$, $d(S_2(t),x')<r$,
contradicting the assumption that $d(x',S_2)=r$.
\QED

\begin{figure}[!ht]
\centering
\includegraphics[width=2.2in]{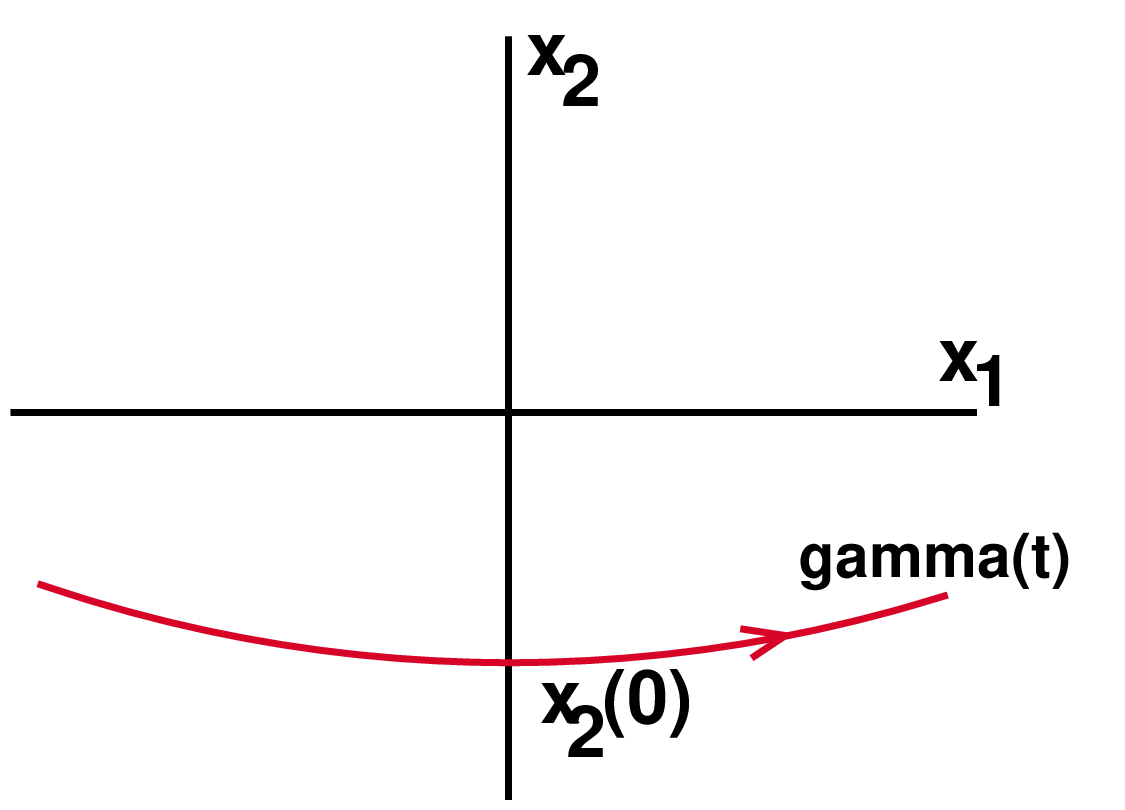}
\caption{\emph{A view of a geodesic $\gamma$ in $K_0$ from `above' (i.e. $x_3>0$). At $t=0$, $\gamma$
passes through the point $(0,x_2(0))$ with velocity $(\dot x_1(0),0)$.}}
\label{fig:funnygeodesic}
\end{figure}

To prove Theorem \ref{thm:funny geodesic2}, we pick a family $\Gamma$ of geodesics in the patch
given by \eqref{eq:coord-patch1} as follows. A geodesic $\gamma(t)$ in $\Gamma$ is determined by initial
condition $\gamma(0)=(0,x_2(0),x_3(0))$, where $x_2(0)$ is not zero but small and $\dot x_1(0)>0$ is of
order unity, while $x_3(0)$ is determined by the fact that $\gamma$ is a curve in the surface $K_0$ (see
Figure \ref{fig:funnygeodesic}).

Since we are interested not in geodesics per se, but in geodesics modulo (non-singular) reparametrization,
we  establish a simple characterization of geodesics in $\Gamma$ that does not depend on the parametrization
(Lemma \ref{lem:curvatureK0}). We then consider the projection $\Pi:K_0\rightarrow K_r$, with $r>0$, which
maps $\gamma$ to a curve $\gamma_r$ in $K_r$. And finally, we prove that $\gamma_r$ is not a
(reparametrization of a) geodesic by showing that it fails the criterion just mentioned.
To do this, we will need to determine the terms of the leading order of magnitude in a fairly involved
expression. We will employ the standard `big-oh' and `small-oh' notation as follows.
We consider curves such as the ones in Figure \ref{fig:funnygeodesic}, and evaluate certain
quantities as these curves cross the $x_1=0$ axis. Thus\footnote{Care should be
exercised with the ``=" sign. It is not reflexive in this context.} using $x$ as shorthand
for $(x_1,x_2)$:
\benn
f(x_1,x_2) &=& Ok \quad \textrm{means} \quad \limsup_{|x|\rightarrow 0} \frac{|f(x)|}{|x|^k} <
\infty \,,\\
\textrm{and  }\quad \quad  f(x_1,x_2) &=& ok \quad \textrm{means} \quad \lim_{|x|\rightarrow 0}
\frac{|f(x)|}{|x|^k}=0 \,.
\eenn

It will be convenient to have a more compact notation. Hence the following definition.

\begin{defn} We define $z_i:=a_ix_i+\partial_i h(x)$, where $x(t)=(x_1(t),x_2(t))$ is the
projection to the $x_1$-$x_2$ plane of the geodesic $\gamma(t)$ in Figure \ref{fig:funnygeodesic}.
\label{def:zi}
\end{defn}

\vskip 0.0in\noindent
We compute the leading orders at $t=0$ of $z_i$, $\dot z_i$, and $\ddot z_i$.
\bsenn
\dot z_i=a_i\dot x_i +d_t\partial_i h \quad \logand \quad \ddot z_i=a_i\ddot x_i +d_t^2\partial_i h \,.
\esenn
We know that $h=o2$ and so $\partial_i h =o1$. Furthermore, at $t=0$, $x_1=0$, and $x_2=O1$. Thus
\bse
z_1=\partial_1 h \quad \logand \quad z_2=a_2x_2 +o1 \,.
\label{eq:zi1}
\ese
Each of these is $O1$ or less. Now,
\bsenn
d_t\partial_i h=\partial_1\partial_ih \;\dot x_1+\partial_2\partial_i h \;\dot x_2 \,.
\esenn
Along the geodesic in the patch, $\dot x_1$ is order unity (or $O1$), and even though $\dot x_2$
may be small, we see that $d_t\partial_2 h=o0$. In fact, we are only interested in evaluating these
quantities at $t=0$ at which point we have $\dot x_2=\ddot x_2 =0$. Putting this together results
at $t=0$ in
\bse
\dot z_1=a_1\dot x_i + o0 \quad \logand \quad \dot z_2 = \partial_1\partial_2 h \;\dot x_1\,,
\label{eq:zi2}
\ese
and so $\dot z_2=o0$. The next derivative, $\ddot z_i$, is a little trickier. The reason is that $d_t^2\partial_2 h$ cannot
be bounded by some order. It may be large, or, depending on $h$, it may be small. To ensure we have
the leading terms of $\ddot z_i$, we have to include both terms and the expression does not simplify.
Setting $t=0$, we will see that $\ddot x_1=O1$ and we know that $\ddot x_1=0$. So at $t=0$,
\bse
\ddot z_1=a_1\ddot x_i + d_t^2\partial_2 h \quad \logand \quad
\ddot z_2 = \partial_1\partial_2 h \;\dot x_1\,.
\label{eq:zi3}
\ese

\vskip 0.2in
\begin{centering}\section{Proof of Theorem \ref{thm:funny geodesic2}}
	\label{chap:funny-geod-example}\end{centering}
\setcounter{figure}{0} \setcounter{equation}{0}

To distinguish the standard inner product in $\R^3$ from
a 2-tuple, we indicate the former by a dot: $x\cdot y$. Also, to avoid cluttering the formulas
with the repetitive occurrence of the argument ``(0)", we will not write it, except when its omission
might lead to misunderstandings.

\begin{lem} Suppose the family of curves $\gamma(t)=(x_1(t),x_2(t), -\frac 12[a_1x_1(t)^2+a_2x_2(t)^2 ]-h)$
in $K_0$ are (a reparametrization of) geodesics with $x_1=0$,
$\dot x_1>0$, $x_2\neq 0$, and $\dot x_2=0$. Then at $t=0$
\bsenn
\lim_{x_2\rightarrow 0}\,\dfrac{\ddot x_2}{\dot x_1^2 x_2}= -a_1a_2 \,.
\esenn
Furthermore, this characterization is independent of the (smooth) parametrization of $\gamma$.
\label{lem:curvatureK0}
\end{lem}

\vskip .1in\noindent
{\bf Proof.} Set $e_i:=\partial_iK_0$, where $K_0$ is given by \eqref{eq:coord-patch1}.
The metric tensor $g_{ij}=e_i\cdot e_j$ and its inverse are given by (see Definition \ref{def:zi})
\bsenn
g=\begin{pmatrix} 1+z_1^2&z_1z_2\\z_1z_2&1+z_2^2\end{pmatrix} \quad \logand \quad
g^{-1}=\Delta^{-1}
\begin{pmatrix} 1+z_2^2&-z_1z_2\\-z_1z_2&1+z_1^2\end{pmatrix} \,,
\esenn
where $\Delta$ is the determinant of $g$. The coefficients of $g^{-1}$ are denoted by $g^{ij}$.
The Christoffel symbols of the second kind are now given by
\bsenn
\Gamma^k_{ij} := \partial_ie_j\cdot \sum_n g^{kn}e_n \,.
\esenn
We have that
\bsenn
\partial_ie_j=(0,0,-\partial_iz_j) \,.
\esenn
So we only need the 3rd component of $\sum_n g^{kn}e_n$. A straightforward computation gives
that these are $-\Delta^{-1}z_k$. This yields
\bsenn
\Gamma^k_{ij} = \Delta^{-1} z_k\;\,\partial_iz_j \,.
\esenn
Employing the rules for order calculation, one checks that this gives an $O1$ term only if $i=j$,
namely $a_ix_i$.
Everything else gives at best $o1$ terms. So $\Gamma^k_{ii} = a_ka_ix_k+o1$, and $\Gamma^k_{ij}=o1$
if $i\neq j$. The geodesic equations are
\bsenn
\ddot x_k+\sum_{i,j}\Gamma^k_{ij}\dot x_i\,\dot x_j =0 \,.
\esenn
So in our case, the equation for $\ddot x_2$ is
\bsenn
\begin{matrix}
\ddot x_2+(a_1a_2x_2+o1)\,\dot x_1^2 + (o1) \dot x_1 \dot x_2+ (a_2^2x_2+o1)\,\dot x_2^2 =0 \,.
\end{matrix}
\esenn
Setting $\dot x_2=0$, proves the first part of the lemma.

To prove that this is invariant under the parametrization $t\rightarrow s(t)$, define
$c(t)=\gamma\circ s$. Set $s(0)=0$. Using $\dot x_2(0)=0$ again, it is trivial to show that at $t=0$
\bsenn
\frac{d_t^2 (x_2(s))}{(d_t x(s))^2\,x_2(s)} = \dfrac{\ddot x_2}{\dot x_1^2 x_2}\,,
\esenn
where we use $d_t$ for $\tfrac{d}{dt}$.
\QED

\begin{lem} Given the surface $K_0$ of \eqref{eq:coord-patch1}, then the constant distance surface
$K_r$ can be parametrized as follows
\bsenn
K_r(u_1,u_2)=\left(u_1,u_2,r-\frac 12(a_{1r}u_1^2+a_{2r}u_2^2) + o2\right)\,,
\esenn
where
\bsenn
a_{ir}=\dfrac{a_i}{1+ra_i} \,.
\esenn
\label{lem:curvatureK0-Kr}
\end{lem}

\begin{figure}[!ht]
	\centering
	\includegraphics[width=1.2in]{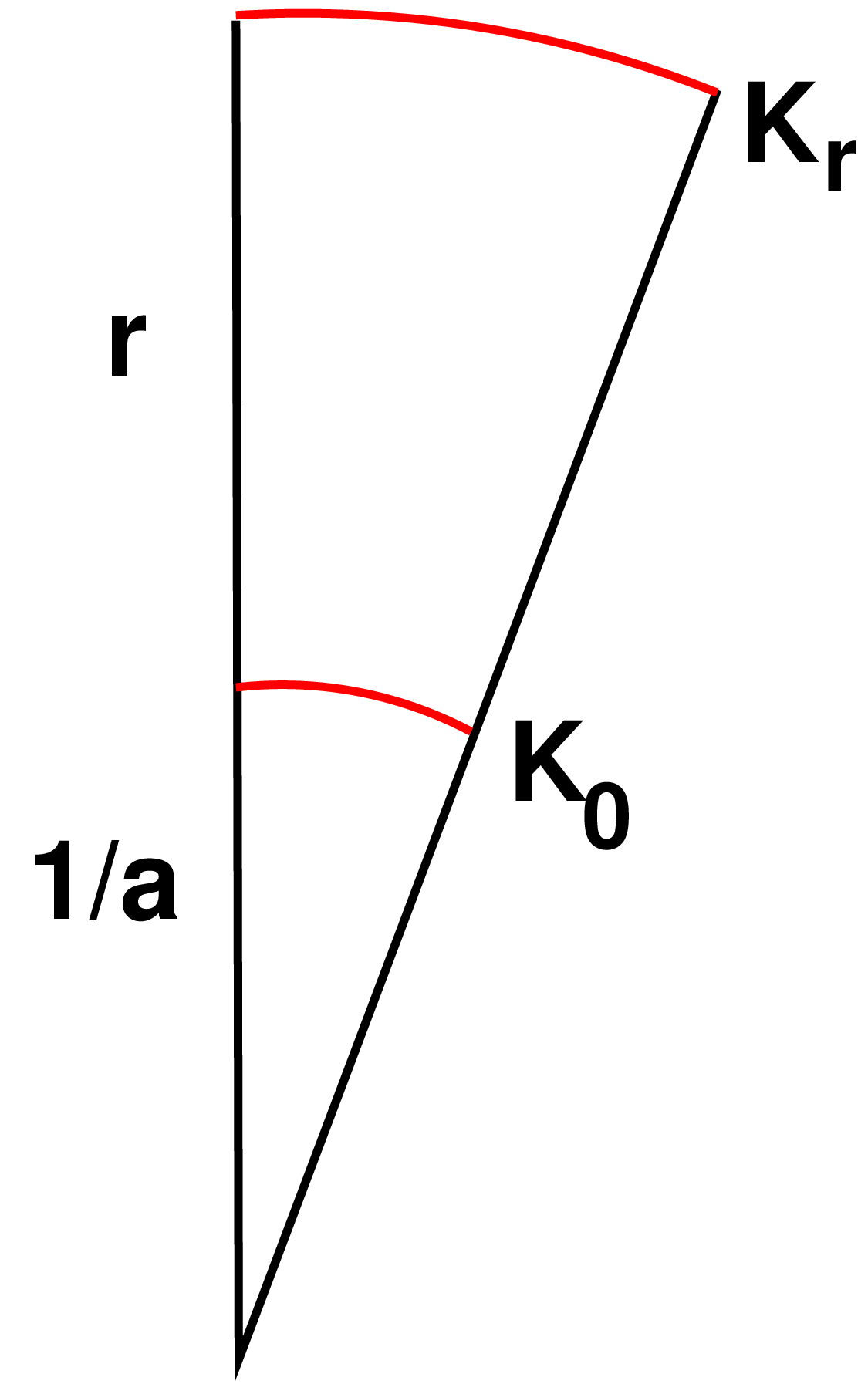}
	\caption{\emph{The radius of curvature in the $x$-direction of $K_0$ at the origin equals $1/a$.
			The orthogonal projection to $K_r$ then gives a radius of curvature of $r+1/a$. The principal
			curvature is the reciprocal of this.}}
	\label{fig:curvature-Kr}
\end{figure}

\vskip-0.1in\noindent
{\bf Proof.} Fix $r>0$. The inverse projection $\Pi_r^{-1}:K_0\rightarrow K_r$ is well-defined, and given by
\bse
K_r(x_1,x_2):=\Pi_r^{-1}(K_0(x_1,x_2))=K_0(x_1,x_2)+r\hat n(x_1,x_2) \,.
\label{eq:inverse-proj}
\ese
We'll call $K_0$, somewhat informally, the `downstairs' surface and $K_r$ is `upstairs'.
We compute, using Definition \ref{def:zi}
\bse
\hat n(x_1,x_2)=\frac{(z_1,z_2, 1)} {\sqrt{1+z_1^2+z_2^2}} \,.
\label{eq:unit-normal}
\ese
So
\benn
K_r(x_1,x_2) &=& \left(x_1+\frac{rz_1}{V}, x_2+\frac{rz_2}V,
- \frac 12(a_1x^2+a_1x_2^2) - h(x_1,x_2) + \frac{r}{V}\right) \,,\\
\where V &=& \sqrt{1+z_1^2+z_2^2}
\eenn
There are no mixed quadratic terms of the form $x_1x_2$ in the expansion of $K_r(x_1,x_2)$.
So if we rewrite this
as $K_r(u_1,u_2)=(u_1,u_2,r+ u_3(u_1,u_2))$, then the $u_1$- and $u_2$-axes of $K_r$ are the axes of
principal curvature at $(u_1,u_2)=(0,0)$. All we
need to do to complete the proof, is a computation of the curvature in the $x_1$-$r$ plane to get
$a_{1r}$. This is done in Figure \ref{fig:curvature-Kr} by employing osculating circles. The
computation is the same in the $x_2$-$r$ plane.
\QED

Part of the difficulty here is that it is pretty clear that if $K_0$ does not have
constant curvature along a geodesic $\gamma(t)$, then the curve traced in $K_r$ by projecting
$\gamma$ will certainly not be a constant speed curve, let alone a constant speed geodesic. It is thus a priori
clear that the projected curve will not satisfy the geodesic equations. What we wish to
establish, however, is whether it can be \emph{reparametrized} as a geodesic. We use Lemma
\ref{lem:curvatureK0} that the images $\gamma_r$ ($r>0$) under the projection are not geodesics.

\begin{lem} The geodesic $\gamma$ depicted in Figure \ref{fig:funnygeodesic} with $\gamma(0)=(0,x_2)$
and $\dot \gamma(0)=(\dot x_1,0)$ projects to a curve
\bsenn
\gamma_r(t)=(u_1(t),u_2(t),u_3(t))
\esenn
in $K_r$ ($r>0$), where at $t=0$, we have
\benn
\dot u_1 &=& (1+ra_1)\dot x_1 +o0\\
u_2 &=& (1+ra_2)x_2+o1\\
\ddot u_2 &=& -(1+ra_1+ra_2)a_1a_2\dot x_1^2x_2+rd_t^2\partial_2 h +o1\,.
\eenn
\label{lem:equations}
\end{lem}

\vskip-0.2in\noindent
{\bf Proof.} We trace a possibly reparametrized geodesic $\gamma(t)$ in $K_0$ satisfying Lemma
\ref{lem:curvatureK0}, and determine the curvature of its projection $\gamma_r$ `upstairs' in $K_r$.
Note that $\gamma_r(t)$ is given by $\gamma(t)+r\hat n(x_1(t),x_2(t))$.
The unit-normal $\hat n$ is given in \eqref{eq:unit-normal}. We use
the rules of evaluating the orders given in Section \ref{chap:outline}.

The $x_1$ and $x_2$ coordinates of $\gamma_r$ will be called $u_1$ and $u_2$ and, noting that $z_i=O1$
(Section \ref{chap:outline}), we get
\bsenn
\begin{matrix}
u_1=x_1+rz_1(1-\frac 12 z_1^2-\frac 12 z_2^2 + O4)\\[0.3cm]
u_2=x_2+rz_2(1-\frac 12 z_1^2-\frac 12 z_2^2 + O4) \,.
\end{matrix}
\esenn
Referring to \eqref{eq:zi1}, this gives for $u_2$ the following:
\bse
u_2=(1+ra_2)x_2+o1 \,.
\label{eq:badcurve1}
\ese
Now, differentiate the $u_i$ with respect to time.
\bsenn
\begin{matrix}
\dot u_1=\dot x_1+r\dot z_1- r\dot z_1(\frac 12 z_1^2+\frac 12 z_2^2+O4)
- r\dot z_1(z_1\dot z_1+z_2\dot z_2+O3)\\[0.3cm]
\dot u_2=\dot x_2+r\dot z_2- r\dot z_2(\frac 12 z_1^2+\frac 12 z_2^2+O4)
- r\dot z_2(z_1\dot z_1+z_2\dot z_2+O3) \,.
\end{matrix}
\esenn
Use \eqref{eq:zi2}, to see that the leading term appearing in $\dot u_1$ is $\dot x_1$  (which is $O0$),
and thus
\bse
\dot u_1=(1+ra_1)\dot x_1 +o0 \,.
\label{eq:badcurve2}
\ese
We need to differentiate $\dot u_2$ one more time with respect to time.
\benn
\ddot u_2 &=& \underbrace{\ddot x_2+\ddot z_2}_A  -\underbrace{r\ddot z_2\;O2}_B -
\underbrace{2r \dot z_2(z_1\dot z_1+z_2\dot z_2+O3)}_C -
\underbrace{rz_2(\dot z_1^2+z_1\ddot z_1+\dot z_2^2+z_2\ddot z_2+O2)}_D\,.
\eenn
To analyze this, we denote the four terms A through D, and look at each individually.
In A, $\ddot z_2$ can be evaluated via \eqref{eq:zi3} and $\ddot x_2$ can be eliminated via Lemma
\ref{lem:curvatureK0}. This gives $-(1+ra_2)a_1a_2\dot x_1^2 x_2 +rd_t^2\partial_2h$ for the term marked
A. Clearly, B is negligible compared to A. In C, $\dot z_2=o0$ as noted before, and the term in
parentheses is $O1$. So all together this term is $o1$ and therefore negligible compared to A. Finally,
in C, we use \eqref{eq:zi1} to establish that $-r\dot z_1^2 z_2$ is $O1$ and, by \eqref{eq:zi1},
all other terms are smaller. This last expression can be simplified using \eqref{eq:zi1} and
\eqref{eq:zi2} to $-ra_1^2a_2\dot x_1^2 x_2+o1$. Collecting terms and adding the relations
\eqref{eq:badcurve1} and \eqref{eq:badcurve2} yields the lemma.       \QED

\vskip-0.2in\noindent
{\bf Proof of Theorem \ref{thm:funny geodesic2}.} On the one hand, if the projected curve $\gamma_r$
is also a geodesic, then it itself must satisfy Lemma \ref{lem:curvatureK0} with the curvatures
given by Lemma \ref{lem:curvatureK0-Kr}. So
\bsenn
\ddot u_2=-\frac{a_1a_2}{(1+ra_1)(1+ra_2)}\dot u_1^2 u_2 \,.
\esenn
Eliminating $\dot u_1$ and $u_2$ in favor of $\dot x_1$ and $x_2$ via Lemma \ref{lem:equations} gives
\bsenn
\ddot u_2=-(1+ra_1)a_1a_2\dot x_1^2 x_2 +o1\,.
\esenn
On the other hand, another equation for $\ddot u_2$ is given by Lemma \ref{lem:equations}.
If we equate the two expressions, we obtain
\bsenn
-(1+ra_1)a_1a_2\dot x_1^2 x_2 +o1= -(1+ra_1+ra_2)a_1a_2\dot x_1^2x_2+rd_t^2\partial_2 h +o1 \,.
\esenn
Upon simplification, this gives
\bse
-  ra_1a_2^2\dot x_1^2x_2 + rd_t^2\partial_2 h = o1 \,.
\label{eq:sphere-or-cylinder}
\ese
This equation has two possible solutions. The first is if the left hand is $o1$ and so $a_1a_2=0$
and $d_t^2\partial_2 h$ is $o1$. From the rules about manipulating the order symbols in Section
\ref{chap:outline}, it follows that then $h=o4$.
The other possibility is if $a_1a_2>0$ and so $d_t^2\partial_2 h=a_1a_2^2\dot x_1^2x_2+o1$.
This happens if and only if $h=\tfrac 14 a_1a_2^2x_1^2x_2^2+g$ and $d_t^2\partial_2 g=o1$.
So $h=\tfrac 14 a_1a_2^2x_1^2x_2^2+o4$.

Now consider the geodesic $\eta$ which is just $\gamma$ rotated by $\pi/2$. Then, by the same reasoning,
if the projection $\eta_r$ in $K_r$ is a geodesic, we must have that $h=\tfrac 14 a_1^2a_2x_1^2x_2^2+o4$.
Since both\footnote{Note that the powers of $a_i$ are distinct.} must hold, we get
$\tfrac 14 a_1a_2^2x_1^2x_2^2=\tfrac 14 a_1^2a_2x_1^2x_2^2$ or $a_1=a_2$.
\QED

\vskip-0.0in\noindent
{\bf Remark.} The two types of solutions of \eqref{eq:sphere-or-cylinder} in this proof do indeed occur.
For if $K_0$ is a plane or a cylinder with radius $1/a$, we get
\bsenn
K(x_1,x_2)=\sqrt{a^{-2}-x_1^2}-a^{-1}=-\tfrac 12 ax_1^2 -\tfrac 18 a^3x_1^4 +O6 \,.
\esenn
In this case, the Gaussian curvature is zero and $h=o4$. On the other hand, for a sphere of radius $1/a$, we have
\bsenn
K(x_1,x_2)=\sqrt{a^{-2}-x_1^2-x_2^2}-a^{-1}=-\frac 12 (ax_1^2+ax_2^2) -\frac 18
(a^3x_1^4+2a^3x_1^2x_2^2+a^3x_2^4) +O6 \,.
\esenn
Here, the principal curvatures are equal and $d_t^2\partial_2 h=a^3\dot x_1^2x_2+o1$.

\vskip 0.2in
\begin{centering}\section{Cylinders Must Be Round}
	\label{chap:round}\end{centering}
\setcounter{figure}{0} \setcounter{equation}{0}

\begin{figure}[!ht]
\centering
\includegraphics[width=2.2in]{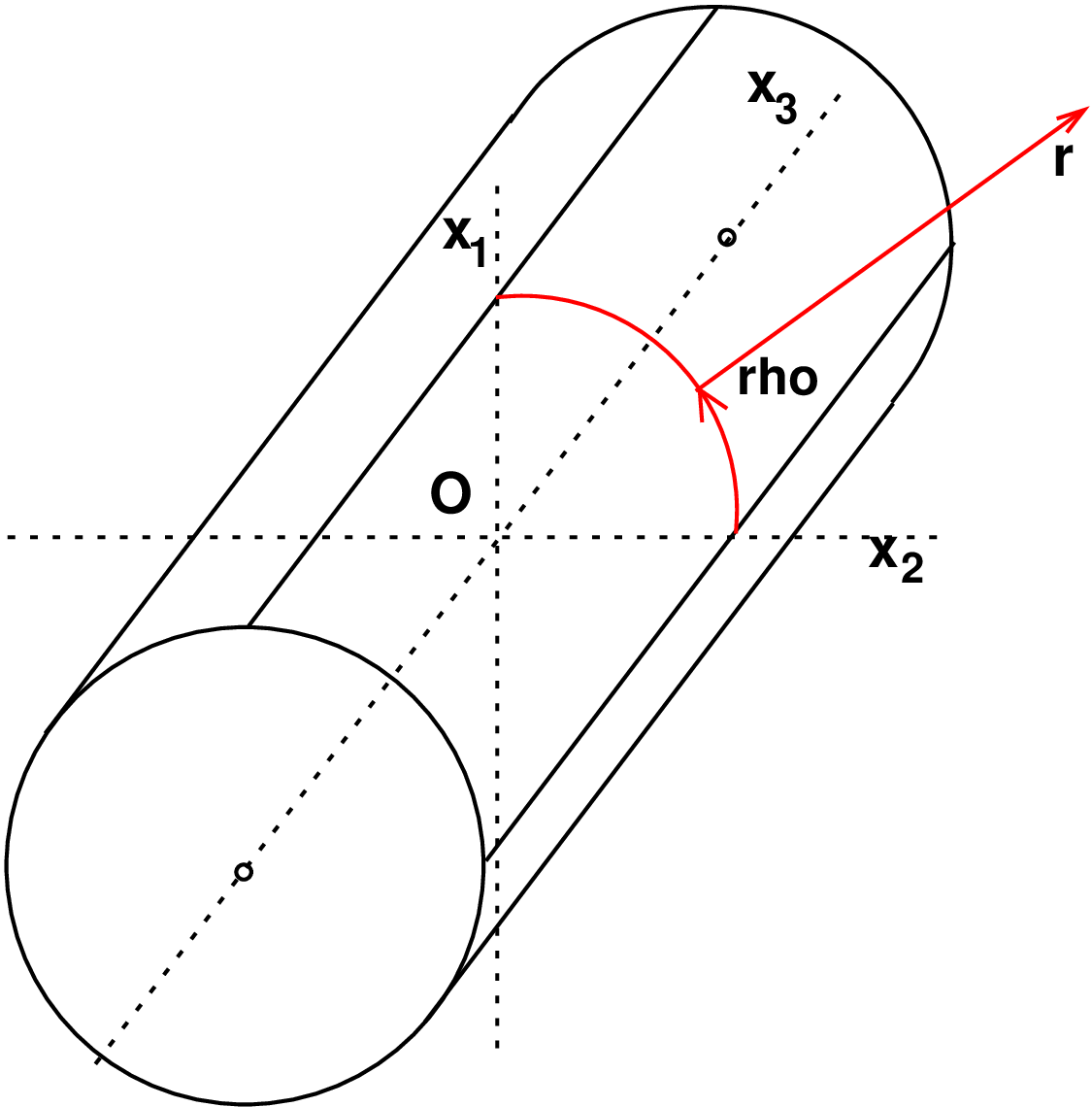}
\caption{\emph{``Polar" coordinates in $\R^3$. }}
\label{fig:polarcoords}
\end{figure}

Now let $K_0$ be a convex, but not necessarily round,  cylinder, invariant under translations along
the $x_3$-axis. Consider ``polar" coordinates $(\rho,x_3,r)$ in $\R^3$  where $\rho$ is the arclength
along the simple, closed curve in the $x_1$-$x_2$ plane that defines $K_0$ as illustrated in Figure \ref{fig:polarcoords}.

\vskip-0.0in\noindent
\begin{prop} With the above assumptions, the projection $\Pi:K_0\rightarrow K_r$
preserves geodesics if and only if the non-zero principal curvature of $K_0$ is constant.
\label{prop:round}	
\end{prop}

\vskip-0.0in\noindent
{\bf Proof.} Since the Gaussian curvature is zero, the map from the cylinder to the plane, given by
$K_0(\rho,x_3)\rightarrow (\rho,x_3)$ is a bijective isometry and so maps geodesics to geodesics.
A geodesic $\gamma$ in $K_0$ is (a) parallel to the $x_3$-axis, or (b) a circle in the $x_1$-$x_2$
plane, or
(c) a curve $\gamma(x_3)=(\rho(x_3),x_3)$. Since $\gamma$ is a geodesic, $\rho(x_3)$ is affine
and has a constant derivative $\frac{d\rho}{dx_3}$. Assume $\gamma$ is a geodesic of type (c).

Now consider the projection $\gamma_r$ of $\gamma$ onto $K_r$. As with $K_0$, we parametrize
$K_r$ by the arclength $\rho_r$ of the defining curve and $x_3$.
It is clear that $\gamma_r$ is a curve $x_3\rightarrow \rho_r(x_3)$. Again, if $\gamma_r$ is a
geodesic, then $\frac{d\rho_r}{dx_3}$ is constant.
Denote the non-zero principal curvature of $K_0$ by $a(\rho)$. A reasoning similar to that
of Lemma \ref{lem:curvatureK0-Kr} gives that arclengths $\rho_r$ and $\rho$ travelled along each
geodesic relate as
\bsenn
d\rho_r=\frac{1/a(\rho)+r}{1/a(\rho)}\,\, d\rho=
(1+a(\rho)r)\,\,d\rho\,.
\esenn
Since the $x_3$ coordinates of $\gamma(t)$ and its projection $\gamma_r(t)$ are the same, we get
\bse
\frac{d\rho_r}{dx_3}=(1+a(\rho)r)\,\,\frac{d\rho}{dx_3}\,.
\label{eq:roundcylinder}
\ese
Thus $\frac{d\rho_r}{dx_3}$ is constant if and only if $a(\rho)$ is constant. \QED

\vskip 0.2in
\begin{centering}\section{A $C^{1,1}$ Counter-example}
	\label{chap:geodesic3}\end{centering}
\setcounter{figure}{0} \setcounter{equation}{0}

\begin{figure}[!ht]
\centering
\includegraphics[width=2.2in]{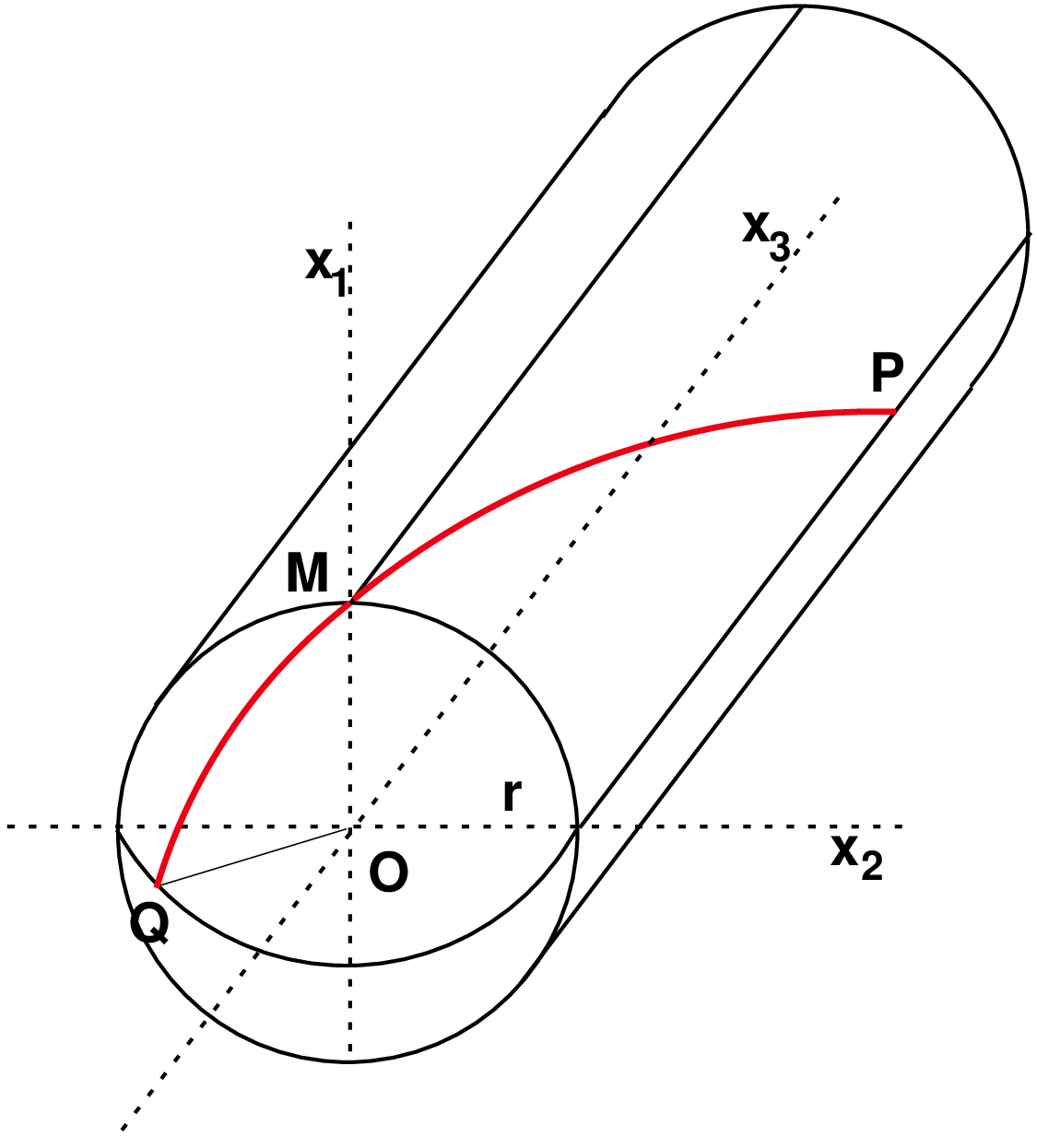}
\caption{\emph{A straight cylinder of radius $r$ topped off by a hemisphere of radius $r$. This surface
is $C^2$, except on the circle where the cylinder and the hemisphere meet. Here $\partial^2/\partial z^2$
has a discontinuity. Since the first derivative changes gradually, this surface is $C^{1,1}$.}}
\label{fig:cylinder-proj}
\end{figure}

We consider the round cylinder `topped off' by a hemisphere both of radius $r$,
which gives a $C^{1,1}$ surface (see Figure \ref{fig:cylinder-proj}). Denote this surface by $S_r$. It
is easy to convince oneself that $S_1$ and $S_r$ ($r>1$) are regular, constant distance surfaces.
Clearly, at every point (except where $C^2$ does not hold) either (i) $a_1a_2=0$ or (ii) $a_1=a_2$.

\vskip 0.1in\noindent
\begin{prop} Let $S_1$, $S_r$, and $\Pi:S_1\rightarrow S_r$ be given as above. Let $\gamma_1$ be the
shortest geodesic connecting $P_1=(0,1,\pi/2)$ and $Q_1=(0,-1/\sqrt{2},-1/\sqrt{2})$. The
projection of $\gamma_1$ by $\Pi$ to $\gamma_r$ (connecting $P_r$ to $Q_r$ in $S_r$ ($r>1$)) is
\emph{not} a local geodesic near the point $M_r$ where $\gamma_r$ intersects the boundary of the
cylinder (see Figure\ref{fig:cylinder-proj}).
\label{prop:funny geodesic3}
\end{prop}

\vskip-0.0in\noindent
{\bf Proof.} The geodesic $\gamma$ connects $P_1$ to $Q_1$, but we do not know where it crosses over
from the cylinder to sphere. So let us call that point $M_1(\theta)$. We have
\bsenn
P_1=(0,1,\pi/2)\;,\quad M_1(\theta)=(\sin \theta, \cos \theta,0) \;,\quad Q_1=\left(0,\tfrac{-1}{\sqrt{2}},\tfrac{-1}{\sqrt{2}}\right) \,.
\esenn
It is easy to see that then the projection $\gamma_r$ of $\gamma$ connects $P_r$ to $Q_r$ via
$M_r(\theta)$ (the same $\theta$), where
\bse
P_r=(0,r,\pi/2)\;,\quad M_r(\theta)=(r\sin \theta, r\cos \theta,0) \;,\quad
Q_r=\left(0,\tfrac{-r}{\sqrt{2}},\tfrac{-r}{\sqrt{2}}\right) \,.
\label{eq:P_r-etcera}
\ese
The projection of $\gamma_r$ consists of two pieces that live on $C^2$ surfaces with either
curvature zero (the cylinder) or the sphere, and so each of these two pieces is a geodesic in $S_r$.

\begin{figure}[!ht]
\centering
\includegraphics[width=4.4in]{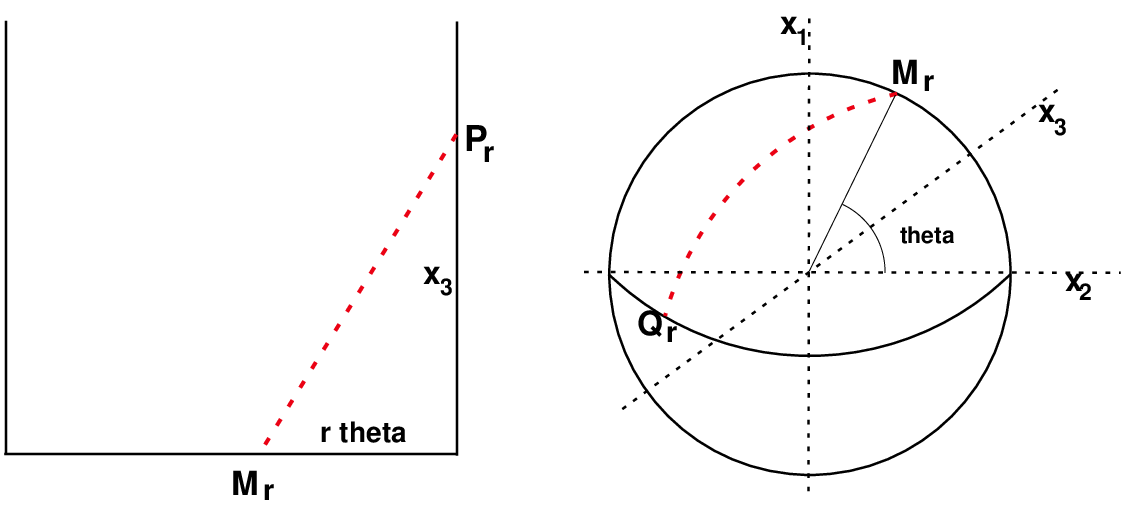}
\caption{\emph{The two geodesic pieces of $\gamma_r$ in $S_r$. To the left, the piece in the
flattened out cylinder. To the right the piece that lies in the hemisphere.}}
\label{fig:piecesofgeodesics}
\end{figure}

The first piece connects $P_r$ to $M_r$, see the left of Figure \ref{fig:piecesofgeodesics}.
In the flattened out cylinder, it is the hypotenuse of the triangle
with sides $\pi/2$ and $r\theta$ and thus has length $\sqrt{\pi^2/4+r^2\theta^2}$. The second
piece lives in the sphere. Its length is $r$ times the angle $\alpha$ between $M_r$ and $Q_r$.
The cosine of $\alpha$ is given by the dot product of the unit vectors parallel to $M_r$ and $Q_r$,
which gives $-\cos \theta/\sqrt{2}$. Thus the length of the second piece equals
$r \arccos \left(\tfrac{-\cos \theta}{\sqrt{2}}\right)$. Therefore, the length of the projected curve
$\gamma_r$ is given by
\bsenn
\ell_\theta(\gamma_r)=\sqrt{\frac{\pi^2}{4}+r^2\theta^2}+r \arccos \left(\frac{-\cos \theta}{\sqrt{2}}\right) \,.
\esenn
We need to minimize this over $\theta\in[0.\pi]$. It is an elementary calculus exercise\footnote{Use
that the derivative of $\arccos q$ equals $-1/\sqrt{1-q^2}$.} to see that this is minimized at
$\theta$ satisfying $\frac{\theta}{\sin \theta}=\frac{\pi/2}{r}$. We know that $\gamma$ is minimizing
in $S_1$. Therefore if we substitute $r=1$, we get $\theta=\pi/2$. That same calculation for $r>1$
implies then that $\gamma_r$ (where $r>1$) is not globally minimizing in $S_r$.

\begin{figure}[!ht]
\centering
\includegraphics[width=2.7in]{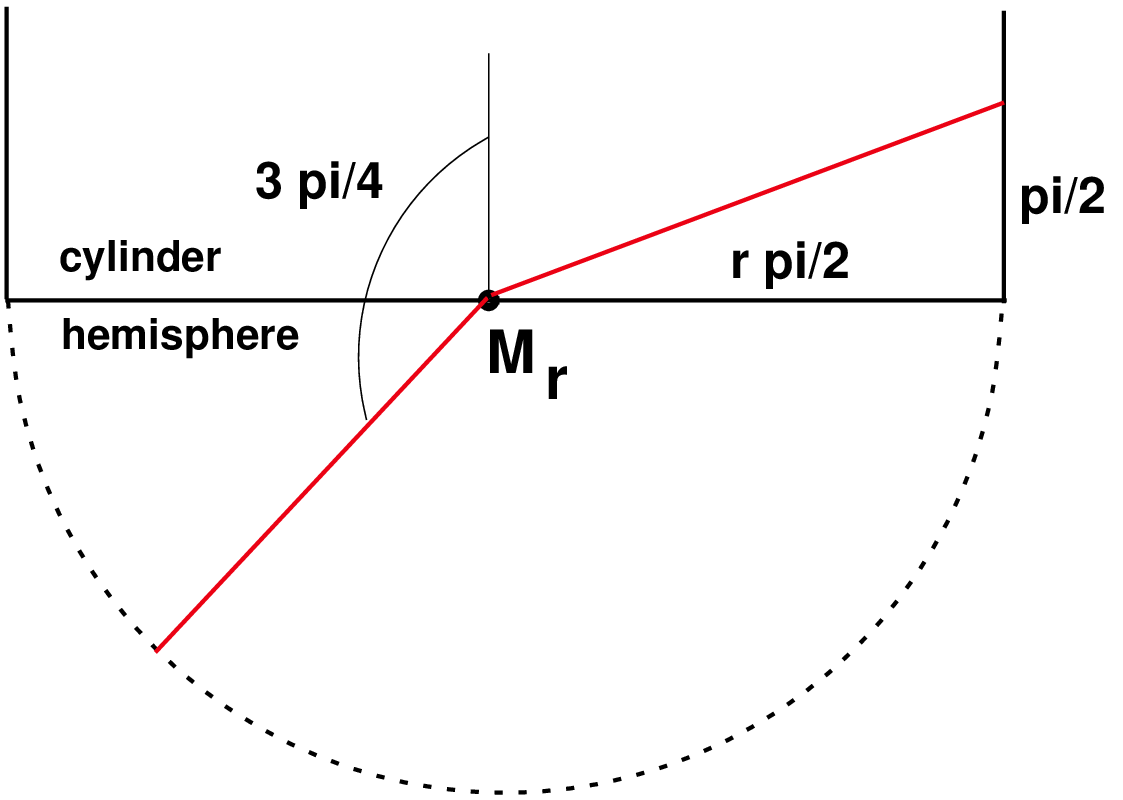}
\caption{\emph{The tangent space $TS_r$ at $M_r$. The red curve corresponds to the lift of $\gamma_r$.
The branch to the left of the base point $M_r$ travels to $Q_r$ and the branch to the right travels
$P_r$ (see \eqref{eq:P_r-etcera}).}}
\label{fig:TS_r}
\end{figure}

Figure \ref{fig:TS_r} is a slightly impressionistic image of the tangent space $TS_r$ at $M_r$.
The slope of $\gamma_r$ restricted to the
lower half plane that projects to the hemisphere equals 1. However the slope restricted the upper
half plane which can be identified with the rolled out cylinder, the slope equals $1/r$.
Thus $\gamma_r$ is not locally minimizing at $M_r$. \QED

\vskip 0.2in
\begin{centering}\section{There are Other Projections that Preserve Geodesics}
	\label{chap:preservegeod}\end{centering}
\setcounter{figure}{0} \setcounter{equation}{0}

In this Section, we find a beautiful example of a family of
projections $\Pi_k:S_k\rightarrow S_0$ such that the surfaces $\left\{S_k\right\}_{k\geq 0}$
foliate the space surrounding the boundary $S_0$ of a convex body in $\R^3$.
It is not known whether this is possible for all such surfaces $S_0$. Our construction is based
on Section \ref{chap:round} and works for (convex) cylinders.

Consider a general, not necessarily round, convex, cylinder $S_0$. It consists of parametrized closed
curve $c(t)$ and lines though that curve, orthogonal to it, as sketched in Figure \ref{fig:polarcoords}.
We can define a $S$ outside $S_0$ by first defining a new curve in $\R^2$:
\bsenn
C(t)=c(t)+r(t)\hat n(t) \,.
\esenn
Here $\hat n(t)$ is the unit normal to $c(t)$ and $r(t)$ is a non-negative distance. Let us denote the
curvature of $c(t)$ by $a(t)$. According to \eqref{eq:roundcylinder}, the projection
$\Pi: S\rightarrow S_0$ between the corresponding cylinders preserves geodesics if
\bsenn
r(t)=\frac{k}{a(t)} \,,
\esenn
where $k$ is a positive constant. The cylinder $S_k$, $k\geq 0$ is given by the lines through $C_k$
orthogonal to plane of $C_k$. Thus the projection $\Pi_k: S_k\rightarrow S_0$ preserves geodesics.
Notice that we have to be careful here, because now the back
and forth projections are not inverses of one another anymore (see Proposition
\ref{prop:projinverses}).

\begin{figure}[!ht]
\centering
\includegraphics[width=5.0in]{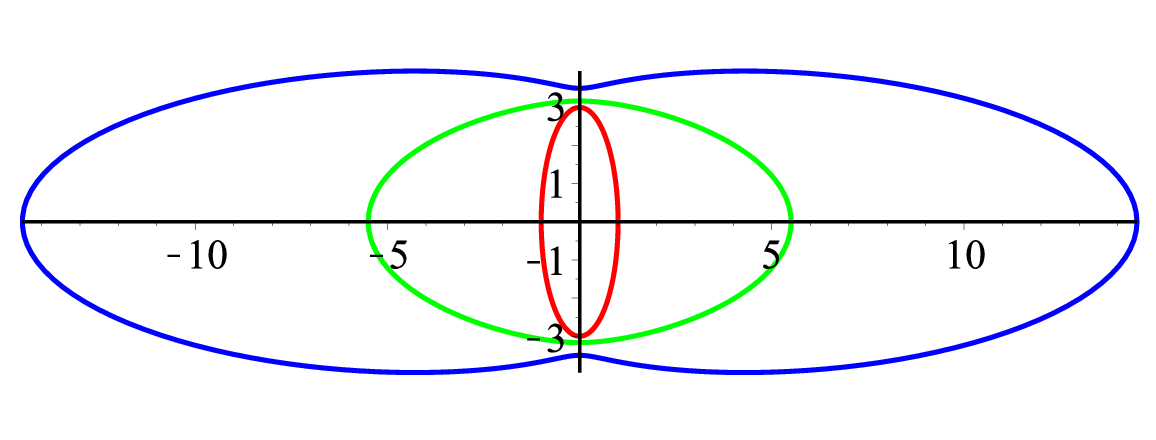}
\caption{\emph{the cylinders orthogonal to the plane of the figure through the blue and
the green curves, have projections to the cylinder orthogonal to the red ellipse that
preserve geodesics.}}
\label{fig:preservinggeodesics}
\end{figure}

We take as an example the ellipse given by
\bsenn
c(t):=(\alpha \cos(t), \beta \sin(t)) \quad \logand \quad C_k(t)=c(t)+\frac{k}{a(t)}\;\hat n(t)\,,
\esenn
where $k$ is a non-negative constant. Standard calculations give $C(t)$ explicitly as
\bsenn
C_k(t)=\left(\left(\alpha+\frac k\alpha(\alpha^2\sin(t)^2+\beta^2\cos(t)^2)\right)\cos(t),
\left( \beta+\frac k\beta(\alpha^2\sin(t)^2+\beta^2\cos(t)^2)\right)\sin(t)\right) \,.
\esenn
We used MAPLE in Figure \ref{fig:preservinggeodesics}, to draw the ellipse $c(t)=(\cos(t), 3 \sin(t))$
in red, $C_k(t)$ for $k=0.5$ in green, and for $k=1.5$ in blue. Note that these remarkable curves
lose convexity for large enough $k$. We leave it to the reader to establish that for large $k$,
the projection $\Pi_k': S_0\rightarrow S_k$ is not single-valued and therefore does not preserve
geodesics.

\vspace{\fill}

\end{document}